\documentclass[11pt,reqno,a4paper,twoside]{article}

\usepackage{amsmath,amsthm,amstext,amscd,amssymb,euscript,url}
\usepackage{mathrsfs}
\usepackage{calrsfs}
\usepackage{epsfig}
\usepackage[inline]{enumitem}

\newcommand{\Z}{\mathbb Z}
\newcommand{\R}{\mathbb R}
\newcommand{\N}{\mathbb N}

\newcommand{\E}{\mathbb E}
\newcommand{\Zd}{\mathbb Z^d}

\renewcommand{\phi}{\varphi}

\newcommand{\loc}{\mathcal{L}}

\def\1{{\mathchoice {\rm 1\mskip-4mu l} {\rm 1\mskip-4mu l}
{\rm 1\mskip-4.5mu l} {\rm 1\mskip-5mu l}}}

\newtheorem{theorem}{{\small T}{\scriptsize HEOREM}}[section]
\newtheorem{corollary}{{\bf{\small C}{\scriptsize OROLLARY}}}[section]
\newtheorem{proposition}{{\bf{\small P}{\scriptsize ROPOSITION}}}[section]
\newtheorem{lemma}{{\bf{\small L}{\scriptsize EMMA}}}[section]
\newtheorem{remark}{{\bf{\small R}{\scriptsize EMARK}}}[section]
\newtheorem{definition}{{\bf{\small D}{\scriptsize EFINITION}}}[section]

\renewenvironment{proof}[1]
{\noindent{{\bf{\small{ P}{\scriptsize ROOF}}}.}\hspace{0.1cm} #1} {$\;\qed$\newline}

\newcommand{\beq}{\begin{eqnarray}}
\newcommand{\eeq}{\end{eqnarray}}

\newcommand{\ba}{\begin{align*}}
\newcommand{\ea}{\end{align*}}

\newcommand{\be}{\begin{equation}}
\newcommand{\ee}{\end{equation}}

\newcommand{\bl}{\begin{lemma}}
\newcommand{\el}{\end{lemma}}

\newcommand{\br}{\begin{remark}}
\newcommand{\er}{\end{remark}}

\newcommand{\bt}{\begin{theorem}}
\newcommand{\et}{\end{theorem}}

\newcommand{\bd}{\begin{definition}}
\newcommand{\ed}{\end{definition}}

\newcommand{\bp}{\begin{proposition}}
\newcommand{\ep}{\end{proposition}}

\newcommand{\bc}{\begin{corollary}}
\newcommand{\ec}{\end{corollary}}

\newcommand{\bpr}{\begin{proof}}
\newcommand{\epr}{\end{proof}}

\newcommand{\bi}{\begin{itemize}}
\newcommand{\ei}{\end{itemize}}

\newcommand{\ben}{\begin{enumerate}}
\newcommand{\een}{\end{enumerate}}

\newcommand{\caD}{{\EuScript D}}

\newcommand{\caH}{{\mathcal H}}

\newcommand{\caK}{{\mathscr K}}

\newcommand{\caS}{{\mathcal S}}

\renewcommand{\(}{\left(}
\renewcommand{\)}{\right)}
\newcommand{\nn}{\nonumber}

\newmuskip\pFqmuskip

\newcommand\pFq[6][8]{%
	\begingroup 
	\pFqmuskip=#1mu\relax
	\mathcode`\,=\string"8000
	\begingroup\lccode`\~=`\,
	\lowercase{\endgroup\let~}\pFqcomma
	{}_{#2}F_{#3}{\left[\genfrac..{0pt}{}{#4}{#5};#6\right]}%
	\endgroup
}
\newcommand{\pFqcomma}{\mskip\pFqmuskip}

\begin{document}
\title{{\bf Quantitative Boltzmann Gibbs principles via orthogonal polynomial duality}}
\author{Mario Ayala, Gioia Carinci and Frank Redig
\\
\small{Delft Institute of Applied Mathematics}\\
\small{Delft University of Technology}\\
{\small Mekelweg 4, 2628 CD Delft}
\\
\small{The Netherlands}
}
\maketitle

\begin{abstract}
We study fluctuation fields of orthogonal polynomials in the context of particle systems with duality.
We thereby obtain a systematic orthogonal decomposition of the fluctuation fields of local functions, where the order of every term can be quantified. This implies a quantitative generalization of the Boltzmann Gibbs principle.
In the context of independent random walkers, we complete this program, including also fluctuation fields in non-stationary context (local equilibrium). For other interacting particle systems with duality such as the symmetric exclusion process,  similar results can be obtained, under precise conditions on the $n$ particle dynamics. 
\end{abstract}

\section{Introduction}
The Boltzmann Gibbs principle is an important ingredient in the study of fluctuation fields of interacting particle systems \cite{kipnis2013scaling}. It basically states that on the central limit scale, the fluctuation field of  local functions can be replaced by a constant times the density fluctuation field, or in other words, it can be replaced by its projection on the one dimensional space generated by the density fluctuation field (where projection has to be understood in an appropriate Hilbert space of macroscopic quantities 
\cite{brox1984equilibrium}). The aim of the present paper is to refine and quantify the Boltzmann Gibbs principle
in the context of particle systems with duality,  using fluctuation fields of orthogonal polynomials.
Indeed, it turns out that replacing the fluctuation field of a local function by its projection on the density field corresponds to the projection on the fluctuation fields of orthogonal polynomials of order one.
Therefore, the Boltzmann Gibbs principle easily follows from an estimation of the covariance of fluctuation fields of orthogonal polynomials of order two and higher. In this paper, for independent random walkers we quantify the precise order of these covariances of fluctuation fields of orthogonal (Charlier) polynomials of order $n$ for all $n\in \N$, and therefore we are able to give an orthogonal decomposition of the fluctuation field of any local function, which is a generalization of the Boltzmann Gibbs principle.
Next, still in the context of independent random walkers, we are able to extend this result in a non-equilibrium setting, using the fact that product of Poisson measures are preserved under this dynamics, i.e., a strong form of propagation of local equilibrium holds in that context.
The basic ingredients of our approach are duality with orthogonal polynomials combined with precise estimates (of local limit type) of the $n$ particle dynamics. Therefore, the results immediately apply in the context of the stationary symmetric exclusion process, and more generally for particle systems where these
precise estimates (of local limit type) of the $n$ particle dynamics can be obtained (e.g. via the log-Sobolev inequality \cite{landim1998decay}).
Next we consider the orthogonal polynomial fluctuation fields themselves and prove that they converge in the sense of generalized processes, i.e., as a random space-time distribution. 
The rest of our paper is organized as follows:
In section 2 we formally introduce our system of random walkers, and the basic concepts and properties needed for the development of this paper. In section 3, on the context of stationarity, we start by introducing our results for the simplest non-trivial example of second order  and move to a generalization first to higher orders and in a next stage to more general functions. We present in section 4 an extension of  these last results to a non-equilibrium setting. Finally in section 5 we show how under additional assumptions our results can be extended to other interacting particle systems.

\section{Basic notions}
\subsection{Independent Random Walkers}
We consider a system of Independent Random Walkers (IRW), an interacting particle system where particles randomly hop on the lattice $\Zd$ without interaction and with no restrictions on the number of particles per site.
Configurations are denoted by $\eta, \xi, \zeta$ and are elements of $\Omega=\N^{\Zd}$ (where
$\N$ denotes the natural numbers including zero).
We denote by $\eta_x$ the number of particles at $x$ in the configuration $\eta\in\Omega$.
The generator working on local functions $f:\Omega\to\R$ is of the type
\be\label{gen}
\loc f(\eta) = \sum_{i,j} p(i,j)\eta_i (f(\eta^{ij})-f(\eta))
\ee
where $\eta^{ij}$ denotes the configuration obtained from $\eta$ by removing a particle from $i$ and putting it at $j$. Additionally, we assume that $p(i,j)$ is a translation invariant, symmetric, irreducible Markov transition function on $\Zd$, i.e.,
\ben
\item $p(i,j)=p(j,i)= p(0,j-i)$.
\item $\sum_{j \in \Zd} p(i,j)=1$
\item There exists $R>0$ such that $p(i,j)=0$ for $|i-j|>R$.
\item For all $x,y\in\Zd$ there exist $i_1=x,\ldots,i_n=y$ such
that $\prod_{k=1}^n p(i_k,i_{k+1})>0$.
\een
For the associated Markov process on $\Omega$, we use the notation $\{\eta(t):t\geq 0\}$, i.e.,
$\eta_x(t)$ denotes the number of particles at time $t$ at location $x \in \Zd$. \\

It is well known that these particle systems have a one parameter family of homogeneous (w.r.t. translations) reversible and ergodic product measures $\nu_\rho, \rho>0$ with Poisson marginals
\[
\nu_\rho (n)= \frac{\rho^n}{n!} e^{-\rho}
\]
This family is indexed by the density of particles, i.e.,
\[
\int\eta_0 d\nu_\rho= \rho
\]

\br
Notice that for these systems the initial configuration has to be chosen in a subset of configurations such that the process $\{\eta(t):t\geq 0\}$ is well-defined. A possible such subset is the set of tempered configurations. This is the set of configurations $\eta$ such that there exist $C, \beta \in \R$ that satisfy $ |\eta (x)| \leq C |x|^\beta$ for all $x\in \R^{d}$. 
We denote this set (with slight abuse of notation) still by $\Omega$, because
we will always start the process from such configurations, and this set has $\nu_\rho$ measure $1$ for all $\rho$. Since
we are working mostly in $L^2(\nu_\rho)$ spaces, this is not a restriction.
\er

\subsection{Orthogonal polynomial self-duality}\label{Dualsec}
The self-duality of the process we introduced and which we need in the sequel is as follows.
We denote by $\Omega_f$ the set of configurations with a finite number of particles (we denote by $\|\xi\|=\sum_x \xi_x$ this number of particles), and
the self-duality function will then be a function $D:\Omega_f\times \Omega \to\R$
such that the following properties hold.
\ben
\item Self-duality:
\be\label{dual1}
\E_\eta \big[D(\xi,\eta_t)\big]=\E_\xi \big[D(\xi_t, \eta)\big]
\ee
for all $\xi\in \Omega_f, \eta\in \Omega$ (where we remind that $\eta\in \Omega$ is always
chosen such that the process $\{\eta(t): t\geq 0\}$ is well-defined when starting from $\eta$).
\item Factorized polynomials:
\be\label{dual2}
D(\xi,\eta)=\prod_{i\in\Zd} d(\xi_i,\eta_i)
\ee
where $d(0,n)=1$, and $d(k,\cdot)$ is a polynomial of degree $k$.
\item Orthogonality:
\be\label{dual3}
\int D(\xi, \eta) D(\xi', \eta) d\nu_\rho(\eta) = \delta_{\xi,\xi'} a(\xi) 
\ee
where $a(\xi)= \|D(\xi, \cdot)\|_{L^2(\nu_\rho)}^2$
\een
Notice that these functions will depend on the parameter $\rho$,  but we suppress this dependence in order not to overload notation.\\

The duality functions which, for Independent Random Walkers, satisfy properties (\ref{dual1}),(\ref{dual2}) and (\ref{dual3})  are known in the literature as Charlier polynomials. These polynomials can be expressed in terms of hypergeometric functions as follows:
\[
d(k,n)= \pFq{2}{0}{-k,-n}{-}{-\frac{1}{\rho}}
\]
the single site duality functions $d \( k, n \)$ satisfy the three terms recurrence relation
\be\label{3terms}
d(k+1,n) =  d(k,n) - \frac{n}{\rho} \; d(k,n-1)
\ee
additionaly to this recurrence relation, at least two more relations can be found. 

\br\label{factrem}  To avoid minor confusions please notice that in \cite{franceschini2017stochastic} a relation between "classical" and new orthogonal duality polynomials is given. Where with classical polynomials we mean
\be\label{classicalpolys}
d(k,n) = \frac{n!}{(n-k)!} 
\ee
and the way they relate is given by
\be\label{chiaform}
D(\xi,\eta ) = \prod_{x \in \Zd}  \sum_{j=0}^{\xi_x} {{\xi_x}\choose{j} } (-\rho)^{\xi_x-j} \frac{\eta_x!}{(\eta_x-j)!} 
\ee
However expression (\ref{chiaform}) differs  by a factor  $-\rho^{|\xi|}$ from the traditional form of the Charlier polynomials found in the literature: 
\be\label{tildeDual}
\tilde{D}(\xi,\eta ) = \prod_{x \in \Zd}  \sum_{j=0}^{\xi_x} {{\xi_x}\choose{j} } (-\rho)^{-j} \frac{\eta_x!}{(\eta_x-j)!} 
\ee
The factor $-\rho^{||\xi||}$  is however invariant under the dynamics of our process that conserves the total number of particles $||\xi(t)||$, and hence its addition preserves the duality property.
Duality fuction \eqref{tildeDual}  is presicely the one that satisfies the relation given in (\ref{3terms}) when starting with $d(0,n)=1$.
\er 
For more details on orthogonal duality and a proof of self-duality with respect to this function we refer to \cite{franceschini2017stochastic} and \cite{redig2017duality}. In those papers a more complete study is provided, which includes the case of other processes such as exclusion and inclusion, among others.\\

\noindent
We denote by $p_t(\xi,\xi')$  the transition probability to go from the  configuration $\xi$ to $\xi'$ in  time $t$.
A key ingredient for  our proof of the Boltzmann Gibbs principle and its extensions is the following elementary consequence of duality with  orthogonal duality functions.
\bl\label{covalemma}
Let $\xi,\xi'\in \Omega_f$, then
\be\label{cova}
\int \E_\eta (D(\xi,\eta_t)) D(\xi',\eta) d\nu_\rho(\eta)= p_t(\xi,\xi') a(\xi')
\ee
\el

\bpr
We use self-duality to compute
\beq
\int \E_\eta [D(\xi,\eta_t)] D(\xi',\eta) d\nu_\rho(\eta) &=& \int \E_\xi [D(\xi_t,\eta)] D(\xi',\eta) d\nu_\rho(\eta) \nn \\
&=& \sum_{\zeta} p_t(\xi,\zeta) \int  D(\zeta,\eta) D(\xi',\eta) d\nu_\rho(\eta) \nn \\
&=& p_t(\xi,\xi') a(\xi')\nn
\eeq
that proves the result.
\epr

\br
Notice that \eqref{cova} in particular implies that if $\eta_0$ is initially distributed according to $\nu_\rho$ then
\be
\text{Cov }_{\nu_\rho} \left( D(\xi,\eta_t) D(\xi',\eta) \right) \geq 0
\ee
i.e.  duality orthogonal polynomials are positively correlated.
\er

\noindent
Lemma \ref{covalemma} provides a big simplification since it allows to transfer most of the uncertainty of our process to the transition kernel $p_t(\xi,\xi')$ of two configurations in $\Omega_f$.  Here $\{\xi(t), t\ge 0\}$ is  a Markov process with countable state space, conserving only $\|\xi(t)\|$ in the course of time, and then easier to treat.  In the Appendix  we provide an estimate of this kernel by means of the local limit theorem. 

\subsection{Fluctuation fields}

Let  $\caS(\R^d)$ be the set of Schwarz functions on $\mathbb R^d$, and denote by $\caS'(\R^d)$
the corresponding distributions space.
Moreover we denote by $\tau_x$ the spatial shift, i.e., $\tau_x (\eta)_y= \eta_{y+x}$,.
Fix $\phi \in \caS(\R^d)$ and let $f:\Omega\to\R$ be  a local function,  we define its fluctuation field on scale $N$ as
\be\label{flucfi}
\mathsf {X_N}(f,\eta;\phi):=a_N(f)\sum_{x\in\Zd} \phi(\tfrac{x}{N})(\tau_x f(\eta)-\psi_f(\rho)) 
\ee
where 
\be
\psi_f(\rho):=\int f d\nu_\rho, \qquad \tau_x f(\eta):= f(\tau_x\eta)
\ee
 and  $a_N(\cdot)$ is a suitable normalization constant depending on $f$.
The field $\mathsf {X_N}(f,\eta;\cdot)$ is a Schwarz-distribution associated to the configuration $\eta$. 
An important case is the density fluctuation field, where we chose $f(\eta)=\eta_0$, $a_N(f)={N^{-d/2}}$.
\vskip.2cm
\noindent
The time-dependent fluctuation field at scale $N$ is then defined as
\be\label{flucfit}
 {\bf \mathsf {X_N}}(f,t;\phi)= \mathsf {X_N}(f,\eta(N^2t);\phi)
\ee
the diffusive rescaling anticipates the natural macroscopic time-scale in this symmetric process, which
has the linear heat equation as hydrodynamic limit.  $\{ {\bf \mathsf {X_N}}(f,t;\cdot), \; t\ge 0\}$ is then a Schwarz-distribution valued stochastic process.

\subsection{Boltzmann-Gibbs principle}
The Boltzmann-Gibbs principle makes rigorous the idea that the density fluctuation field is the fundamental fluctuation field, because the
density is the only (non-trivial) conserved quantity in the process under consideration. This means that one can replace, in first approximation,
the fluctuation field of a function $f$ by its ``projection on the density field''.
For a local function $f$ this projection is the fluctuation field of  the function $P_1(f):=\psi_f'(\rho) (\eta_0-\rho)$, where $\psi_f(\rho) = \int f d \nu_\rho$.

\vskip.2cm
\noindent
The standard statement of the Boltzmann Gibbs principle is given in the following theorem.
\bt\label{TEO:BG}
For all $f$ local, and $\phi\in\caS(\R^d)$ and for all $T>0$
\be\label{bomgi}
\lim_{N\to\infty}\E_{\nu_\rho} \left[\frac 1 {N^{d/2}} \int_0^T\( {\bf \mathsf {X_N}}(f,t;\phi)- {\bf \mathsf {X_N}}(P_1(f),t;\phi)\) dt\right]=0
\ee
\et
\noindent
We refer to \cite{kipnis2013scaling} for the proof of Theorem \label{TEO:BG} and for a comprehensive discussion of the result that is valid in a more general context and not only for the process  considered in the present paper.

\subsection{Fluctuation fields of orthogonal polynomials}

For $n\in\N$ we denote by $\caH_n$ the (real) Hilbert spaces generated by the polynomials $D(\xi,\cdot)$ with degree at most $n$, i.e. $||\xi||\le n$. We have of course the inclusion $\caH_0=\R\subset \caH_1\subset\caH_2\subset\ldots$ and the union of the spaces $\caH_n$ is dense in $L^2(\nu_\rho)$. Moreover, for every $f\in L^2(\nu_\rho)$ its projection on $\caH_n$ is given by
\be\label{proj}
f_n= \sum_{\xi\in\Omega_f: \|\xi\|\leq n} \langle f, D(\xi,\cdot)\rangle \frac{D(\xi,\cdot)}{a(\xi)}
\ee
where $\langle \cdot ,\cdot \rangle$ denotes the $L^2 (\nu_\rho)$ inner product.\\ 

\noindent
The aim of what follows is to show that the Boltzmann Gibbs principle is an instance of a more general statement concerning the fluctuation behavior of functions which are orthogonal to   $\caH_n$ for some $n\in \mathbb N$. This is (in some sense to be explained below) the case for the function $f-P_1(f)$. 
\vskip.2cm
\noindent
For $\xi\in\Omega_f$, $\phi\in\caS(\R^d)$ we define the n-th order polynomial  fluctuation field as
\beq
X_N(\xi,\eta,\phi) &&:=\sum_{x\in \mathbb Z^d}\phi\(\tfrac x N\) \, D( \xi, \tau_x \eta) \nn \\
&&=\sum_{x\in \mathbb Z^d}\phi\(\tfrac x N\) \, D(\tau_x \xi,  \eta)
\eeq

\vskip1cm

\section{Stationary case}
\subsection{Second order polynomial field}
We start with the simplest non-trivial example for independent random walkers started from a product measure with homogeneous Poisson marginals.
To illustrate our point let us start with a simple computation, which contains all the important ingredients
of the more general Theorem \ref{thmirw} below. 
Consider the field
\be
\label{1fimacvar}
X_N^{(2)}(\eta;\phi):=X_N(2\delta_0,\eta,\phi)=\sum_{x \in \Zd}\phi\(\tfrac{x}{N}\)\,  D(2\delta_x, \eta)
\ee
The notation $X_N^{(2)}$ suggests that this is in some sense the ''second order" polynomial field.
In the orthogonal polynomial language, this is the field of the second order Charlier polynomial:
\be
D(2\delta_x, \eta)=\eta_x(\eta_x-1)-2\rho(\eta_x-\rho)-\rho^2
\ee
recall from earlier that
\[
a ( 2 \delta_0) = \int(D(2\delta_x,\eta))^2 d\nu_\rho(\eta)
\]
then we have the following.
\bp\label{2ficovProp}
The second order polynomial field $X_N^{(2)}(\eta;\phi)$ is such that
\ben
\item For $t>0$ we have
\be\label{2ficov}
\E_{\nu_\rho}\left[X_N^{(2)}(\eta(t);\phi)\, X_N^{(2)}(\eta(0);\phi)\right]=a ( 2 \delta_0) \,\sum_{x,y \in \Zd}\phi(\tfrac{x}{N})\phi(\tfrac{y}{N}) (p_t(x,y))^2
\ee
\item As a consequence, for $t>0$ we have
\be\label{10fimacvar}
\lim_{N\to\infty}\E_{\nu_\rho}\left[X_N^{(2)}(\eta(N^ 2t);\phi) X_N^{(2)}(\eta(0);\phi)\right] = \frac{d\cdot a ( 2 \delta_0) }{(2\pi t)^d}\int_{\R^{2d}} e^{-\frac{d|x-y|^2}{t}} \phi(x)\phi(y) dx dy
\ee
\een
\ep
\bpr
The first statement follows from self-duality and  Lemma \ref{covalemma}.
For the second statement we use that $\phi$ has compact support, call this support $S$, and define 
\be\label{Mdef}
M := \max \{ d(x,y) : x, y \in S\}
\ee
it follows from Theorem \ref{LCLTCRW} that there exists $c= c(M)$ such that 
\[
\sup_{x:|x|\le M N \sqrt{t}}  p^{RW}_{N^2t} (x)  \le \bar p_{N^2t}(x) \(1+ \frac c { N\sqrt t}\)
\]
with $ \bar p_{t}(\cdot)$ as defined in \eqref{pbar}. Then from \eqref{2ficov} it follows that
\beq
&&\E_{\nu_\rho}\left[X_N^{(2)}(\eta(t);\phi)X_N^{(2)}(\eta(0);\phi)\right] \nn \\
&&=
a ( 2 \delta_0) \sum_{x,y \in S}\phi(\tfrac{x}{N})\phi(\tfrac{y}{N})\,\bar p_{N^2t}(x) \bar p_{N^2t}(y) \(1+ \frac c { N\sqrt t}\)^2 \nn \\
&&=
a ( 2 \delta_0) \cdot \frac {d} {(2\pi t)^{d}}\cdot  \frac 1 {N^{2d}}\sum_{x,y \in S}\phi(\tfrac{x}{N})\phi(\tfrac{y}{N}) e^{-\frac{d(z-y)^2}{tN^2}}  \(1+ \frac c { N\sqrt t}\)^2 \nn
\eeq
and letting $N \to \infty$ we obtain the r.h.s. of \eqref{10fimacvar}.
\epr
\noindent
In the current context the Boltzmann Gibbs principle for the fluctuation field of the function $ f = \eta_0 ( \eta_0 -1)$ is a consequence of Proposition \ref{2ficovProp}. We make this statement more transparent with the following corollary 
\bc\label{BGcor}
The field $X_N^{(2)}(\eta(N^ 2 t);\phi) $ is such that for all $T>0$ and for all $N$ big enough
\be\label{Bint}
\frac1{N^d}\int_0^T\int_0^T\E_{\nu_\rho}\left[X_N^{(2)}(\eta(N^ 2 t);\phi) X_N^{(2)}(\eta(N^2 s);\phi)\right]\, ds \, dt\leq C(T) N^{-\frac{2d}{ 2 + d}}
\ee
More precisely, \eqref{10fimacvar} gives a better estimate of the order of the covariance of the fluctuation field in the diffusive time-scale as $N\to\infty$.
\ec

\bpr
Given the fact that the RHS of \eqref{10fimacvar} has an indetermination at $t=0$. Hence we derive the following estimate for the integrand in \eqref{Bint}
\beq
&& \frac1{N^d}\E_{\nu_\rho}\left[ X_N^{(2)}(\eta(N^ 2 t);\phi) X_N^{(2)}(\eta(N^2 s);\phi)\right] \nn \\
&&= K_\rho\, \frac1{N^d} \sum_{x \in \Zd}\phi(\tfrac{x}{N}) p_{N^2(t-s)}(x,y) \sum_{y \in \Zd} \phi(\tfrac{y}{N}) p_{N^2(t-s)}(x,y) \nn \\
&&\leq K_{\rho} p_{N^2(t-s)}(0,0)  \| \phi \|_1 \E_x \phi(\tfrac{X_t}{N})  \nn \\
&&\leq K_{\rho} p_{N^2(t-s)}(0,0)  \| \phi \|_1   \| \phi \|_{\infty} \nn
\eeq
at this point we could have concluded \eqref{Bint} by naively estimating $ p_{N^2(t-s)}(0,0) $ by one. Nevertheless our aim is to provide a more quantitative statement. Hence, we distinguished the cases $ |t-s| \geq \epsilon_N$ and $ |t-s| < \epsilon_N$ where $\epsilon_N$ is to be optimized. By the LCLT
\be
p_{N^2(t-s)}(0,0) \leq \frac {d} {(2\pi N^2(t-s))^{d/2}}
\ee
then
\be
  p_{N^2(t-s)}(0,0) \leq
    \begin{cases}
      \frac {d} {N^d \epsilon_N^{d/2}}, & \text{if}\ |t-s| \geq \epsilon_N \\
      1  & \text{if}\ |t-s| < \epsilon_N
    \end{cases}
\ee
Hence the integral is bounded by
\beq
&&\int_0^T\int_0^T \frac1{N^d}\E_{\nu_\rho}\left[ X_N^{(2)}(\eta(N^ 2 t);\phi) X_N^{(2)}(\eta(N^2 s);\phi)\right]\, ds \, dt \nn \\
&&\leq K_{\rho}  \| \phi \|_1   \| \phi \|_{\infty} \frac{T^2}{2} \left[ \frac {d} {N^d \epsilon_N^{d/2}} + d \epsilon_N  \right]
\eeq
Assume $\epsilon_N $ is of the form $ N^{-\alpha}$, optimality then comes from solving for $\alpha$ 
\[
N^{-\alpha} = N^{-d} N^{d/2 \alpha}
\]
after elementary computations we find $\alpha = \frac{2d}{d+2}$. Which in fact not only shows that the Boltzmann-Gibbs principle holds, but also provides us with a better estimate of the order of convergence. 
\epr

Back to the second order polynomial fluctuation fields, and for the sake of transparency, we make explicit the dependency on the "coordinate points" $x_1, x_2 $ and redefine the fields in terms of the orthogonal duality polynomials as follows:
\be
\label{2ndfield}
X_N^{(2)}(x_1, x_2,\eta;\phi):=\sum_{x \in \Zd}\phi(\tfrac{x}{N}) D(\delta_{x_1+x}+\delta_{x_2+x} , \eta)
\ee
Notice then, that in Proposition \ref{2ficovProp} we treated  for $x_1=x_2=0$. It is necessary then to verify that Proposition \ref{2ficovProp} is not only result of this particular choice we made,  consider then for $x_1\not=x_2$ the field
\be
X_N^{(2),\not=}(x_1, x_2,\eta,\phi)= \sum_{x \in \Zd} \phi(\tfrac{x}{N})(\eta_{x+x_1}-\rho)(\eta_{x+x_2}-\rho) 
\ee
where the upper index $\not=$ refers to the fact that $x_1\not= x_2$. We then have the following analogous of Proposition \ref{2ficovProp}.
\bp\label{2diffficovprop}
The second order  polynomial fluctuation field $X_N^{(2),\not=}(x_1, x_2,\eta;\phi)$ is such that
\ben
\item For $t>0$ we have
\beq\label{2ficov}
&&\E_{\nu_\rho}(X_N^{(2),\not=}(x_1, x_2,\eta(t);\phi)X_N^{(2),\not=}(x_1, x_2,\eta(0);\phi)) \nonumber \\
&=& a( \delta_{x_1} + \delta_{x_2})\sum_{x,y \in \Zd}\phi(\tfrac{x}{N})\phi(\tfrac{y}{N}) p_t(x+x_1,x+x_2;y+x_1,y+x_2) \nonumber\\
&+& a( \delta_{x_1} + \delta_{x_2})\sum_{x,y \in \Zd}\phi(\tfrac{x}{N})\phi(\tfrac{y}{N})p_t(x+x_1,x+x_2;y+x_2,y+x_1) 
\eeq
\item As a consequence, for $t>0$ we have
\beq\label{2fimacvar}
&& \lim_{N\to\infty} \E_{\nu_\rho}(X_N^{(2),\not=}(x_1, x_2,\eta(N^ 2t);\phi) X_N^{(2),\not=}(x_1, x_2,\eta(0);\phi)) \nn \\
&&=\frac{2 a( \delta_{x_1} + \delta_{x_2}) d}{(2\pi t)^d}\int_{\R^{2d}} e^{-\frac{d|x-y|^2}{t}} \phi(x)\phi(y) dx dy
\eeq
\een
\ep
\bpr
The argument for the first statement is similar to the one in the proof of Proposition \ref{2ficovProp}, the difference is that now
\[
D(\delta_{x+x_1} +\delta_{x+x_2},\eta)= (\eta_{x+x_1}-\rho)(\eta_{x+x_2}-\rho) 
\]
is the product of two first order Charlier polynomials, which by the assumption of factorized polynomials allows us to proceed in the same way than before. Furthermore, in this case we have
\beq
&&p_t(\delta_{x+x_1}+ \delta_{x+x_2},\delta_{y+x_1}+ \delta_{y+x_2}  ) \nn \\
&&= p_t(x+x_1,x+x_2;y+x_1,y+x_2) + p_t(x+x_1,x+x_2;y+x_2,y+x_1) \nn \\
\eeq
which is the source of the second term in \eqref{2ficov}. In the second statement is necessary to verify that $x_1$ and $x_2$ do not play a role in the leading order
\beq\label{exprimacvar}
&&\E_{\nu_\rho}(X_N^{(2),\not=}(x_1, x_2,\eta(N^ 2t);\phi) X_N^{(2),\not=}(x_1, x_2,\eta(0);\phi)) \nonumber \\
&=&a( \delta_{x_1} + \delta_{x_2})\sum_{x,y \in \Zd}\phi(\tfrac{x}{N})\phi(\tfrac{y}{N}) p_{N^2t}(x+x_1,x+x_2;y+x_1,y+x_2) \nonumber \\
&+&a( \delta_{x_1} + \delta_{x_2})\sum_{x,y \in \Zd}\phi(\tfrac{x}{N})\phi(\tfrac{y}{N}) p_{N^2t}(x+x_1,x+x_2;y+x_2,y+x_1)
\eeq
The first term in the RHS of \eqref{exprimacvar} can be treated in the same way than before. For the second term, we just have to notice 
\[
|x+x_1-y-x_2|^2 +|x+x_2-y-x_1|^2  = 2|x-y|^2 +2|x_1-x_2|^2
\]
and proceed in the same way. 
\epr

Now we show how to  generalize  this result and  discuss the case of higher order fields.

\subsection{Higher order fields}

Let $k\in \mathbb N$ and denote by $\mathbf x\in \mathbb Z^{kd}$ the coordinates vector $\mathbf x:=(x_1,\ldots,x_k)$, with $x_i\in {\mathbb Z}^d$, $i=1,\ldots, k$. We denote by $\xi(\mathbf x)$ the configuration associated to $\mathbf x$, i.e. $\xi_x(\mathbf x)=\sum_{i=1}^k\mathbf 1_{x=x_i}$. We define  $||\mathbf x||:=||\xi(\mathbf x)||=k$. Here $x_i$ is the position of the $i$-th particle, where particles are labeled in such a way that the dynamics is symmetric. For a more extensive explanation of the labeled dynamics we refer the reader to \cite{demasi2006mathematical}.
We  denote by $\hat \tau_z$, $z\in \mathbb Z^d$ the shift operator acting on the coordinate representation:
\be
\hat {\tau}_z \mathbf x= (z+x_1, \ldots, z+x_k), \qquad \text{and then} \qquad   \tau_z\xi = \xi(\hat \tau_z \mathbf x)
\ee
Because of the translation invariance of the dynamics we have that 
\be\label{qua}
p_t(\xi(\hat \tau_y \mathbf x), \xi(\hat \tau_z \mathbf x))=p_t(\xi(\mathbf x), \xi(\hat \tau_{z-y} \mathbf x))
\ee
With an abuse of notation, we keep denoting by $p_t(\mathbf x, \mathbf y)$ the transition probability of the labeled particles in the coordinate representation. 
\br The relation between the transition probabilities in the coordinate and in the configuration representations is given by 
\be\label{transprob}
p_t(\xi(\mathbf x),\xi(\mathbf y))=\sum_{\mathbf x': \\ \xi(\mathbf x')=\xi(\mathbf y)} p_t(\mathbf x, \mathbf x')
\ee
\er
Notice that it is presicely from  relation \eqref{transprob} that a factor of 2 appears in Proposition \ref{2diffficovprop} and not in Proposition \ref{2ficovProp}. We can expect that in this general setting the difference among cases will become more cumbersome. To avoid any further notational difficulties we introduce the following: \\

Let $\mathcal P_k$ be the set of permutations of $\{1,\ldots, k\}$, for $\sigma,\sigma' \in \mathcal P_k   
$ we define the following equivalence relation:
\be
\sigma \sim \sigma'  \quad \text{mod} \quad \mathbf x \qquad \text{iff} \quad x_{\sigma(i)}=x_{\sigma'(i)} \quad \forall i \in \{1,\ldots,k\}
\ee
and define $\mathcal P_k(\mathbf x):=\mathcal P_k/\sim_{\mathbf x}$. Then we have
\be
|\mathcal P_k(\mathbf x)|=\frac{k!}{\prod_{i\in \mathbb Z^d}\xi_i(\mathbf x)!}
\ee
For each $\sigma\in \mathcal P_k(\mathbf x)$ we define the new coordinate vector $\mathbf x^{(\sigma)}$ such that
\be
\mathbf x^{(\sigma)}_i= x_{\sigma(i)}
\ee
thus we can write
\be
p_t(\xi(\mathbf x), \xi(\hat \tau_{z} \mathbf x))=\sum_{\mathbf x': \\ \xi(\mathbf x')=\xi(\hat\tau_z \mathbf x)} p_t(\mathbf x, \mathbf x')= \sum_{\sigma\in \mathcal P_k({\mathbf x})}p_t(\mathbf x, \hat \tau_z \mathbf x^{(\sigma)})
\ee
With a slight abuse of notation we denote by
\be
X_N(\mathbf x,\eta,\phi):= \sum_{z\in \mathbb Z^d} \phi \( \frac z N\) D(\hat \tau_z \mathbf x,\eta), 
\ee
 define the $k$-th order fluctuation field associated to the $k$-particles configuration $\mathbf x$. Then  we have
\bt\label{thmirw}
Let $k:=||\mathbf x||$, then the $k$-th order fluctuation field $X_N(\mathbf x,\eta,\phi)$ is such that
\ben
\item
For all $t>0$
\beq\label{kficov}
&&\mathbb E_{\nu_\rho}\left[X_N(\mathbf x,\eta(t),\phi)X_N(\mathbf x,\eta(0),\phi)\right] \nn\\
&&= a(\xi(\mathbf x))\sum_{\sigma\in \mathcal P_k({\mathbf x})} \sum_{y,z\in \mathbb Z^d} \phi\(\frac y N\)\phi\(\frac z N\)p_t( \mathbf x, \hat \tau_{z-y} \mathbf x^{(\sigma)})
\eeq
\item
As a consequence, for $t>0$
\beq\label{kfimaccov}
&&\lim_{N\to\infty} N^{d(k-2)}\mathbb E_{\nu_\rho}\left[X_N(\mathbf x,\eta(N^2t),\phi)X_N(\mathbf x,\eta(0),\phi)\right] \nn\\
&=&
|\mathcal P_k(\mathbf x)| a(\xi(\mathbf x))  \frac{ d^{k/2}}{(2\pi t)^{dk/2}}\int_{\R^{2d}} e^{-kd|z-y|^2/2t} \phi(z)\phi(y) dz dy
\eeq
\een
\et

\bpr
The first statement of the theorem is a direct application of Lemma \ref{covalemma}  and the fact that the function $a(\cdot)$ is translation invariant, i.e. $a(\xi(\hat \tau_z \mathbf x))=a(\xi(\mathbf x))$, for all $z \in \mathbb Z^d$.
\beq\label{qui}
&&\mathbb E_{\nu_\rho}\left[X_N(\mathbf x,\eta(t),\phi)X_N(\mathbf x,\eta(0),\phi)\right] \nn\\
&&= a(\xi(\mathbf x))\sum_{y,z\in \mathbb Z^d} \phi\(\frac y N\)\phi\(\frac z N\)p_t(\xi(\hat \tau_y \mathbf x), \xi(\hat \tau_z \mathbf x))
\eeq
Then, from \eqref{qua} and \eqref{qui} it follows that
\beq
&&\mathbb E_{\nu_\rho}\left[X_N(\mathbf x,\eta(t),\phi)X_N(\mathbf x,\eta(0),\phi)\right]\nn\\
&&= a(\xi(\mathbf x))\sum_{\sigma\in \mathcal P_k({\mathbf x})} \sum_{y,z\in \mathbb Z^d} \phi\(\frac y N\)\phi\(\frac z N\)p_t( \mathbf x, \hat \tau_{z-y} \mathbf x^{(\sigma)})
\eeq
For the second stament observe that from translation invariance we have
\be\label{transobs}
p^{\text{IRW}}_{N^2t}( \mathbf x, \hat \tau_{z-y} \mathbf x)=\(p_{N^2t}^{RW}(z-y)\)^k
\ee
Define $B_{M,N}:=\{x\in \mathbb Z^d: |x|\le N M\}$, then, since $\phi$ has a finite support we have that there exists $M\ge 0$ such that, for
\beq
&&\sum_{y,z\in \mathbb Z^d} \phi\(\frac y N\)\phi\(\frac z N\)p^{\text{IRW}}_{N^2t}( \mathbf x, \hat \tau_{z-y} \mathbf x)\nn\\
&& =\sum_{y,z\in B_{M,N}} \phi\(\frac y N\)\phi\(\frac z N\) \(p_{N^2t}^{RW}(z-y)\)^k\nn\\
&&=\(\frac {\sqrt d} {(2\pi t)^{d/2}}\)^k\(1+\frac{c}{N\sqrt t}\)^k\frac 1 {N^{kd}}\sum_{y,z\in  B_{M,N}} \phi\(\frac y N\)\phi\(\frac z N\) \, e^{-\frac{kd|\frac z N-\frac y N|^2}{2t}}\nn
\eeq
for a suitable $c=c(M)$, the last inequality coming from Theorem \ref{LCLTCRW}. We have 
\beq
&&\lim_{N\to \infty} \frac 1 {N^{2d}}\sum_{y,z\in \mathbb Z^d} \phi\(\frac y N\)\phi\(\frac z N\) e^{-\frac{kd|\frac z N-\frac y N|^2}{2t}}= \int_{\mathbb R^{2d}} \phi\(y\)\phi\(z\) \, e^{-\frac{kd|z -y|^2}{2t}} dxdz \nn
\eeq
\epr

\subsubsection{Quantitative Boltzmann-Gibbs principle}
 
On the same spirit than Corollary \ref{BGcor} we can now state a refined quantitative version of the Boltzmann-Gibbs principle for higher order fields. 

\bt\label{BGRk}
The field $X_N^{(k)}(\eta(N^ 2 t);\phi) $ is such that for all $T>0$ there exists $C(T)$  such that  for all $N$ big enough
\be\label{Bintk}
\frac1{N^d}\int_0^T\int_0^T\E_{\nu_\rho}\left[ X_N(\mathbf x,\eta(N^2 t),\phi)X_N(\mathbf x,\eta(N^2 s),\phi)  \right]\, ds \, dt \leq C(T) N^{-\frac{2(k-1) d}{ 2 +(k-1) d}}
\ee
\et

\bpr
Analogously to the case of two particles ( see the proof of Corollary \ref{BGcor}), and using observation \eqref{transobs} we first obtain the following estimate
\beq
&& \frac1{N^d} \E_{\nu_\rho}\left[ X_N(\mathbf x,\eta(N^2 t),\phi)X_N(\mathbf x,\eta(N^2 s),\phi)  \right] \nn \\
&&\leq  \(p_{N^2(t-s)}^{RW}(0)\)^{k-1} |\mathcal P_k(\mathbf x)| a(\xi(\mathbf x))  \| \phi \|_1   \| \phi \|_{\infty} 
\eeq
again, by the LCLT 
\be
 \(p_{N^2(t-s)}^{RW}(0)\)^{k-1}  \leq
    \begin{cases}
      \frac {d} {N^{(k-1)d} \epsilon_N^{(k-1)d/2}}, & \text{if}\ |t-s| \geq \epsilon_N \\
       1, & \text{otherwise}
    \end{cases}
\ee
allowing us to bound the integral
\beq
&&\int_0^T\int_0^T \frac1{N^d}\E_{\nu_\rho}\left[ X_N^{(2)}(\eta(N^ 2 t);\phi) X_N^{(2)}(\eta(N^2 s);\phi)\right]\, ds \, dt \nn \\
&&\leq |\mathcal P_k(\mathbf x)| a(\xi(\mathbf x))  \| \phi \|_1   \| \phi \|_{\infty} \frac{T^2}{2} \left[ \frac {d} {N^{(k-1)d} \epsilon_N^{(k-1)d/2}} + d \epsilon_N  \right]
\eeq
the same anzats, $\epsilon_N = N^{-\alpha}$, results on the optimal value
\be
\alpha = \frac{2(k-1) d}{ 2 +(k-1) d}
\ee
\epr

\subsection{Fluctuation Fields of projections on $\caH_N$}

We can further generalize part (2) of Theorem \ref{thmirw} to a wider class of functions $f$. In this section we make such a generalization for a particular subset of $L^2(\nu_\rho)$. For $f\in L^2(\nu_\rho)$  we can use the fact that the union of the spaces $\caH_n$ is dense in $L^2(\nu_\rho)$  to express $f$ as follows
\be\label{denexp}
f (\eta) = \sum_{\substack{ n \geq 0 \\ \xi\in \Omega_f: \|\xi\|= n}} C_{n,\xi} D(\xi,\eta)
\ee
for the rest of this section we restrict ourselves to the set of functions $f\in L^2(\nu_\rho)$ satisfying the following condition
\be\label{expancondi}
\sum_{  \xi, \xi^\prime \in \Omega_f: \|\xi\|=  \|\xi^\prime\|} | C_{n,\xi} C_{n,\xi^\prime} | a(\xi^\prime) < \infty
\ee
In particular all linear combinations of orthogonal duality polynomials satisfy \eqref{expancondi}. 

\bt\label{gralirw}
Let $f$ be a function such that the condition \eqref{expancondi} is satisfied, and as before let $f_{k-1}$ denote the projection of $f$ on $\caH_{k-1}$, then the field
\[
\mathsf {X_N}(f-f_{k-1},\eta;\phi)=\sum_{x\in \Zd} (\tau_x f (\eta)- \tau_x f_{k-1}(\eta)) \phi\(\frac x N\)
\]
satisfies
\[
\E_{\nu_\rho}\left[\mathsf{X_N}(f-f_{k-1},\eta;\phi)\mathsf{X_N}(f-f_{k-1},\eta(N^2 t);\phi)\right] = O(N^{-d(k-2)})
\]
\et

\bpr
After some simplifications due to orthogonality the field reads 
\[
\mathsf{X_N}(f-f_{k-1},\eta;\phi)= \sum_{x \in \Zd} \phi\(\frac x N\)  \sum_{\substack{ n \geq k \\ \xi\in\Omega_f: \|\xi\|= n}} C_{n,\xi}  \tau_x D(  \xi,\eta) 
\]
We then compute
\beq 
&&\E_{\nu_\rho}\left[\mathsf{X_N}(f-f_{k-1},\eta;\phi)\mathsf{X_N}(f-f_{k-1},\eta(N^2 t);\phi)\right] \nonumber \\
&=&
\sum_{x, y} \phi\(\frac x N\) \phi\(\frac y N\) \sum_{\substack{ n \geq k \\ \xi\in \Omega_f: \|\xi\|= n}} \sum_{\substack{ l \geq k \\ \xi^\prime \in \Omega_f: \|\xi^\prime \|= l}} C_{n,\xi} C_{l,\xi^\prime} \int \tau_x D(\xi,\eta) E_\eta \left[ \tau_y D(\xi^\prime,\eta(N^2 t)) \right] d \nu_\rho (\eta) \nonumber \\
&=&
\sum_{x, y} \phi\(\frac x N\) \phi\(\frac y N\)\sum_{\substack{ n \geq k \\ \xi\in \Omega_f: \|\xi\|= n}} \sum_{\substack{ l \geq k \\ \xi^\prime \in \Omega_f: \|\xi^\prime \|= l}} C_{n,\xi} C_{l,\xi^\prime} \int \tau_x D(\xi,\eta) E_\eta \left[ \tau_y D(\xi^\prime,\eta(N^2 t)) \right] d \nu_\rho (\eta) \nonumber \\
&=&
\sum_{x, y} \phi\(\frac x N\) \phi\(\frac y N\) \sum_{\substack{ n \geq k \\ \xi\in \Omega_f: \|\xi\|= n \\ \xi^\prime \in \Omega_f: \|\xi^\prime\|= n}} C_{n,\xi} C_{n,\xi^\prime} a(\xi^\prime)   p_{N^2 t}(\tau_y \xi^\prime, \tau_x \xi)
\eeq
from the LCLT we can also obtain that
\[
p_{N^2 t}(\tau_y \xi, \tau_x \xi^\prime) = \mathcal{O}(N^{-d \|\xi\|})
\]
this, allows us to bound our expression of interest
\beq 
&&N^{d(k-2)} \E_{\nu_\rho}\left[\mathsf{X_N}(f-f_{k-1},\eta;\phi)\mathsf{X_N}(f-f_{k-1},\eta(N^2 t);\phi)\right] \nn \\
&\leq&
N^{d(k-2)} \sum_{x, y} \phi\(\frac x N\) \phi\(\frac y N\)   \sum_{\substack{ n \geq k \\ \xi\in \Omega_f: \|\xi\|= n \\ \xi^\prime \in \Omega_f: \|\xi^\prime\|= n}}  \frac{M}{N^{dn}} | C_{n,\xi} C_{n,\xi^\prime}| a(\xi^\prime) \nn \\
&=&
\left( \frac{1}{N^{2d}}\sum_{x, y} \phi\(\frac x N\) \phi\(\frac y N\) \right) \sum_{\substack{ n \geq k \\ \xi\in \Omega_f: \|\xi\|= n \\ \xi^\prime \in \Omega_f: \|\xi^\prime\|= n}} \frac{M}{N^{d(n-k)}} | C_{n,\xi} C_{n,\xi^\prime}| a(\xi^\prime) \nn \\
\eeq
At this point we need to show that the last summation does not play a role in the leading order. But this comes from the fact  that  $f$ satisfies condition \eqref{expancondi}. 
\epr
Analogously to Theorem \ref{BGRk} we provide a quantitative version of the Boltzmann-Gibbs principle for the current setting. 

\bt\label{ProjBGR}
The field $\mathsf {X_N}(f-f_{k-1},\eta;\phi)$ is such that for all $T>0$ there exists $C(T)$ such that for all $N$ big enough
\be\label{Bintk}
\frac1{N^d}\int_0^T\int_0^T\E_{\nu_\rho}\left[ \mathsf {X_N}(f-f_{k-1},\eta(N^2 t);\phi) \mathsf {X_N}(f-f_{k-1},\eta (N^2 s);\phi)  \right]\, ds \, dt \leq C(T) N^{-\frac{2(k-1) d}{ 2 +(k-1) d}}
\ee
\et

\section{Non-stationary fluctuation fields}
\subsection{Second order  fields}
Let us now start independent walkers from a product measure of  non-homogeneous Poisson, with weakly varying density profile
i.e., from $\otimes_{x\in\Zd} \nu_{\rho(\tfrac{x}{N})}$.
We denote by $\caD_\rho$ the orthogonal polynomials, i.e., 
\[
\caD_\rho(\xi,\eta)= \prod_{i} \caD_{\rho(i)}(\xi_i,\eta_i)
\]
where $\caD_{\rho(i)}$ denote the orthogonal polynomials w.r.t. Poisson with parameter $\rho(i)$.
We now are interested in the fields
\be
X_N(\xi,\rho,\phi,t):= \sum_{x \in \mathbb Z^d} \phi(\tfrac{x}{N})  \caD_{\rho_{tN^2}} (\xi,\eta (N^2 t)) 
\ee
 then the second order field is
\be
X_N^{(2)}(\rho,\phi,t):=X_N(2\delta_0,\rho,\phi,t)= \sum_{x}\phi(\tfrac{x}{N}) \caD_{\rho_{tN^2}} (2\delta_x, \eta (N^2 t))
\ee
with respect to previous notation please notice the additional dependence on the parameter $\rho$, where $\rho_t(x)= \E_x\left[\rho(X_t)\right]$ and $X_t$ denotes the continuous-time random walk.\\

We want to prove that the covariance of $X_N^{(2)}(\rho,\phi,t)$ and $X_N^{(2)}(\rho,\phi,s)$ is of order 1, as $N\to\infty$, exactly as in the stationary case. For this we start with the following result:
\bl\label{orthodual}
Let $\nu_\rho:=\otimes_{x\in\Zd} \nu_{\rho(x)} $ be a product of non-homogeneous Poisson measures, then
we have
\be
\int \E_\eta \left[\caD_{\rho_t(x)} (2\delta_x, \eta(t))\right] \caD_{\rho(y)} (2\delta_y,\eta) \; d\nu_\rho (\eta)= k_2(y)\,  p_t(x,y)^2
\ee
where 
\[
k_2 (y)=\int\left( \caD_{\rho(y)} (2\delta_y,\eta)\right)^2 d\nu_\rho (\eta)
\]
\el

\bpr
Note that
\[
\caD_{\rho_t(x)} (2\delta_x, \eta_t)= \eta_x(t)(\eta_x(t)-1)-2\rho_t(x)(\eta_x(t)-\rho_t(x)) -\rho_t(x)^2
\]
hence
\beq
&&\E_\eta \left[\caD_{\rho_t(x)} (2\delta_x, \eta_t)\right]
\nonumber\\
&=&
\E_\eta \left[\eta_x(t)(\eta_x(t)-1)\right] -2\rho_t(x)\E_\eta\left[\eta_x(t)-\rho_t(x)\right]-\rho_t(x)^2
\eeq
We now state the following:\\
{\bf Claim 1:}\\
\[
\int \E_\eta\left[\eta_x(t)-\rho_t(x)\right]\caD_{\rho(y)} (2\delta_y,\eta) \, d\nu_\rho (\eta)=0
\]
Indeed, by duality, $\E_\eta\left[\eta_x(t)-\rho_t(x)\right]=\sum_z p_t(x,z) (\eta_z-\rho(z))$
and $(\eta_z-\rho(z))$ is in $L^2(\nu_\rho (\eta))$ always orthogonal to 
$\caD_{\rho(y)} (2\delta_y,\eta)$ because for $z\neq y$ both $(\eta_z-\rho(z))$ and
$\caD_{\rho(y)} (2\delta_y,\eta)$ have expectation zero and when $z=y$ because
it is the inner product of the first order and second order orthogonal polynomials, which
is zero.
So we only have to work out the expectation $\E_\eta \left[\eta_x(t)(\eta_x(t)-1)\right]$
which by duality equals
\[
\sum_{u} p_t(x,u)^2 \eta_u(\eta_u-1) + 2\sum_{u\not= v}p_t(x,u) p_t(x,v) \eta_u \eta_v
\]
{\bf Claim 2:}\\
For all $u$
\[
\int \eta_u\caD_{\rho(y)} (2\delta_y,\eta) d\nu_\rho (\eta)=0
\]
Indeed, for $u\not= y$ this is true because of the  product character of the measure and
the fact that $\caD_{\rho(y)} (2\delta_y,\eta)$  has zero expectation, and  for
$u=y$ $\eta_y=\eta_y-\rho(y) + \rho(y)$ which is the sum of the first orthogonal polynomial and
a constant, which is in $L^2(\nu_\rho (\eta))$ orthogonal to
$\caD_{\rho(y)} (2\delta_y,\eta)$.

Finally, we remark that for all $u\not=y$ 
\[
\int \eta_u(\eta_u-1)\caD_{\rho(y)} (2\delta_y,\eta) d\nu_\rho (\eta)=0
\]
because of the  product character of the measure and
the fact that $\caD_{\rho(y)} (2\delta_y,\eta)$  has zero expectation.
Finally, 
\[
\int \eta_y(\eta_y-1)\caD_{\rho(y)} (2\delta_y,\eta) d\nu_\rho (\eta)= \int (\caD_{\rho(y)} (2\delta_y,\eta))^2 d\nu_\rho (\eta)
\]
because adding first order terms in $\eta_y$  does not change the inner product with $\caD_{\rho(y)} (2\delta_y,\eta)$.
\epr
\\
As a consequence of Lemma \ref{orthnb odual} and using that a product of Poisson measures is reproduced at later times, we compute 
\beq
&&\lim_{N\to\infty}\E_{\nu_\rho} \left[X_N^{(2)}(\rho,\phi,t) X_N^{(2)}(\rho, \phi, s)\right]
\nonumber\\
&=& \lim_{N\to\infty}\E_{\nu_{\rho_{sN^2}}} \left[X_N^{(2)}(\rho,\phi,t-s) X_N^{(2)}(\rho, \phi, 0)\right]
\nonumber\\
&=&
\int \frac{e^{-\frac{(x-y)^2}{t-s}}}{2\pi (t-s)^{d/2}}\; \phi(x)\phi(y)  \, \kappa_2(y) dx dy
\eeq
where
\[
\kappa_2(y)= \lim_{N\to\infty} k_2(Ny)
\]
which exists because the initial Poisson measure has slowly varying density profile.

\subsection{Higher order fields: Non-stationary case}
The aim of this section is to extend the results of the previous example to higher order fields:
\be
X_N(\mathbf x,\rho,\phi,t)= \sum_{x \in \mathbb Z^d} \phi(\tfrac{x}{N})  \caD_{\rho_{tN^2}} (\hat \tau_{\mathbf x}  \xi , \eta (N^2 t)) 
\ee
We start then with a generalization of Lemma \ref{orthodual} to higher orders. As we already stated in Remark \ref{factrem} in the case of independent random walkers, the orthogonal duality polynomials are related to the classical duality polynomials in the following way:
\be
\caD_{\rho}(\xi,\eta ) = \prod_{x \in \Zd}  \sum_{j=0}^{\xi_x} {{\xi_x}\choose{j} } (-\rho(x))^{\xi_x-j} d (j, \eta_x ) 
\ee
where $d(k,n)$ are the classical single site duality polynomials. 
\br 
Notice that due to the non-homogeneity of the product measure, the duality property cannot be any longer guaranteed.
\er
Despite of the previous remark, the special form of the Charlier polynomials allows us to reach the same conclusions than in the stationary case. Let us first make a simple observation:\\
Define $A(\xi, \eta,\rho)$ as the difference between the Charlier and classical polynomials of order $\|\xi\|$, i.e.
\[
A(\xi, \eta,\rho) := \caD_{\rho}(\xi,\eta ) - \prod_{x \in \Zd}  d (\xi_x, \eta_x )
\]
and notice that $A(\xi, \eta,\rho)$ is a polynomial of degree strictly less than $\|\xi\|$ and as a consequence it has an expansion, in temrs of orthogonal polynomials, consisting only on polynomials of order strictly smaller than $\|\xi\|$. Therefore, by orthogonality we have
\[
\int \E_\eta \left[ A(\xi, \eta,\rho)  \right] \caD_{\rho_0} (\xi^\prime,\eta) d\nu_{\rho_0} (\eta)= 0
\]
for any configuration $\xi^\prime$ such that $\|\xi\| \leq \|\xi^\prime\|$. With this observation we are ready to state the following Lemma:
\bl\label{HorthoDual}
Let $\nu_\rho:=\otimes_{x\in\Zd} \nu_{\rho(x)} $ be a product of non-homogeneous Poisson measures, and
let $\rho_t(x)= \E_x\left[\rho(X_t)\right]$, where $X_t$ denotes continuous-time random walk. Then
we have
\be
\int \E_\eta \left[\caD_{\rho_t} (\xi, \eta(t))\right] \caD_{\rho_0} (\xi^\prime,\eta) d\nu_{\rho_0} (\eta)= p_t(\xi,\xi^\prime) a_0(\xi^\prime)
\ee
where $a_t(\xi)= \|\caD_{\rho_t}(\xi, \cdot)\|_{L^2(\nu_{\rho_t})}^2$
\el
\bpr
We simply compute
\beq
&&\int \E_\eta \left[\caD_{\rho_t} (\xi, \eta(t))\right] \caD_{\rho_0} (\xi^\prime,\eta) d\nu_{\rho_0} (\eta) \nn \\
&=& 
\int \E_\eta \left[\prod_{x }   \sum_{j=0}^{\xi_x} {{\xi_x}\choose{j} } (-\rho_t)^{\xi_x-j} d (j, \eta(x,t) ) \right] \caD_{\rho_0} (\xi^\prime,\eta) d\nu_{\rho_0} (\eta) \nn \\
&=& 
\int \E_\eta \left[ \prod_{x}  d (\xi_x, \eta(x,t) ) \right] \caD_{\rho_0} (\xi^\prime,\eta) d\nu_{\rho_0} (\eta) + \int \E_\eta \left[A(\xi, \eta,\rho)  \right] \caD_{\rho_0} (\xi^\prime,\eta) d\nu_{\rho_0} (\eta) \nn \\
&=& 
\int \E_\xi \left[ \prod_{x}  d (\xi(x,t), \eta_x ) \right] \caD_{\rho_0} (\xi^\prime,\eta) d\nu_{\rho_0} (\eta) \nn \\
&=& 
\int \sum_\zeta p_t(\xi, \zeta) \left( \prod_{x }  d (\zeta_x, \eta_x ) + A(\zeta, \eta,\rho) \right) \caD_{\rho_0} (\xi^\prime,\eta) d\nu_{\rho_0} (\eta) \nn \\
&=& 
\int \sum_\zeta p_t(\xi, \zeta) \left( \prod_{x}   \sum_{j=0}^{\zeta_x} {{\zeta_x}\choose{j} } (-\rho_0)^{\zeta_x-j} d (j, \eta_x ) \right)  \caD_{\rho_0} (\xi^\prime,\eta) d\nu_{\rho_0} (\eta) \nn \\
&=& 
\int \sum_\zeta p_t(\xi, \zeta) \caD_{\rho_0} (\zeta,\eta)   \caD_{\rho_0} (\xi^\prime,\eta) d\nu_{\rho_0} (\eta) \nn \\
&=&
p_t(\xi,\xi^\prime) a_0(\xi^\prime)
\eeq
where in the fourth and fifth line we subtracted and added zero respectively by  using the orthogonality of $\caD_{\rho_0} (\xi^\prime,\eta)$ to lower order polynomials in the expansion.
\epr
We now state the non-stationary version of Theorem \ref{thmirw}
\bt
\label{NSthmirw}
Let $\nu_\rho:=\otimes_{x\in\Zd} \nu_{\rho(x)} $ and $\rho_t(x)$ be as before, and let $k:=||\mathbf x||$, then
\ben
\item For all $t>0$
\beq
&&\E_{\nu_\rho}\left[X_N(\mathbf x,\rho,\phi,t)  X_N(\mathbf x ,\rho,\phi,0)\right] \nn\\
&=&
a_0\left(\sum_{i=1}^k\delta_{x_i} \right) \sum_{x,y} \phi(\tfrac{x}{N})\phi(\tfrac{y}{N}) p_t \left(\sum_{i=1}^k\delta_{ x+x_i}; \sum_{i=1}^k \delta_{y+x_i}\right)
\eeq
\item
As a consequence, for $t>s>0$
\beq
&&\lim_{N\to\infty} N^{d(k-2)}\E_{\nu_\rho}\left[X _N(\mathbf x,\rho,\phi,t)  X_N(\mathbf x,\rho,\phi,s) \right]
\nn \\
&=&
\caK(x_1,\ldots, x_k;\rho) \frac{d^{k/2}}{(2\pi (t-s))^{dk/2}}\int_{\R^2} e^{-kd(x-y)^2/2(t-s)} \phi(x)\phi(y) dx dy \nn
\eeq
with $ \xi = \sum_{i=1}^k\delta_{x_i} $ and $\caK(x_1,\ldots, x_k;\rho)$ defined as in the stationary case.
\een
\et
\bpr
Is a consequence of Lemma \ref{HorthoDual} together with the fact that a product of Poisson measure is reproduced at later times. 
\epr

With this last theorem, we have now the ingredients to obtain a quantitative Boltzmann-Gibbs principle 
\bc
For all $T>0$there exists $C(T)$ such that for all $N$ big enough
\be\label{BGLocEq}
\frac1{N^d}\int_0^T\int_0^T\E_{\nu_\rho}\left[ X _N(\mathbf x,\rho,\phi,t) X _N(\mathbf x,\rho,\phi,s)  \right]\, ds \, dt \leq C(T) N^{-\frac{2(k-1) d}{ 2 +(k-1) d}}
\ee
\ec

\bpr
The proof is essentialy the same than in all the previous cases.
\epr

\section{Particle systems with Orthogonal Duality}

In the context of stationarity, the results of this paper are not exclusive for Independent Random Walkers. Hence in this section we extend our results to a wider class of IPS. i.e. to those particle systems that enjoy the existence of orthogonal self-duality and that satisfy an additional condition in the transition kernel. Let then  $\{  \eta_t  \}_{ t \geq 0} $ be an IPS for which there exists an orthogonal self-duality function $D: \Omega_f \times \Omega \to \R $ satisfying all the properties stated in section \ref{Dualsec}. As in the same section, we denote by $p_t(\xi,\xi^\prime)$ the transition probability to go from configuration $\xi$ to $\xi^\prime$ in time $t$. Then, inmmediatly follows the following

\bl\label{covalemmaNONLCLT}
Let $\xi,\xi'\in \Omega_f$, then
\be\label{covaNONLCLT}
\int \E_\eta (D(\xi,\eta_t)) D(\xi',\eta) d\nu_\rho(\eta)= p_t(\xi,\xi') a(\xi')
\ee
\el

\noindent
furthermore, let us assume that for all $ \xi,\xi'\in \Omega_f$, the transition kernel satisfies the following estimate
\be\label{kerassum}
p_t (\xi, \xi^\prime) \leq \frac{C}{(1+t)^{ \| \xi \|d/2 }} 
\ee
This assumption is reasonable, since in \cite{landim1998decay} estimates of this kind were already found for a wide class of interacting particle systems that for example includes generalized exclusion processes. The results of \cite{landim1998decay} are applicable as long as the process satisfies a logarithmic Sobolev inequality for the symmetric part of the generator.
\noindent
As before, for a fix $\mathbf x \in \Z^{dk}$ we define the polynomial fluctuation field
\be
X_N(\mathbf x,\eta,\phi):= \sum_{z\in \mathbb Z^d} \phi \( \frac z N\) D(\hat \tau_z \mathbf x,\eta),
\ee
from assumption \eqref{kerassum} we can also conclude
\bt
For all $T>0$ there exists $C(T)$ such that for all $x \in \Z^{dk}$ and for all $N$ big enough
\be\label{BGNonLCLT}
\frac1{N^d}\int_0^T\int_0^T\E_{\nu_\rho}\left[ X_N(\mathbf x,\eta(N^2 t),\phi)X_N(\mathbf x,\eta(N^2 s),\phi)  \right]\, ds \, dt \leq C(T) N^{-\frac{2(k-1) d}{ 2 +(k-1) d}}
\ee
\et

\section{Appendix}

\subsection{Local limit theorems}

In this section we state and prove a local central limit theorem for Independent Random Walkers in continuous time. The motivation of this section comes from the fact that, despite of being common knowledge, we were not able to find a reference that includes the proof of such a result. However we do have access to many versions of the discrete case. We state now the version included in \cite{lawler2010random}, since we consider is the most suitable to then jump to the continuous time case. Theorem \ref{LCLTDRW} below is a direct consequence of Theorem 2.1.1 in the same reference \cite{lawler2010random}.
\vskip.4cm
\bt[LCLT for Discrete-Time  Random Walk]\label{LCLTDRW}
Let $x\in \mathbb Z^d$ and  $p^{\text{DRW}}_n(\cdot)$ be the probability distribution of a discrete-time random walk in $\mathbb Z^d$, then, for any fixed $M\ge 0$ there exists $c=c(M)$ such that 
\be
\sup_{|x|\le M\sqrt n}\bigg|\frac{p^{\text{DRW}}_n(x)}{\bar p_n(x)}-1 \bigg| \le \frac {c} {n}
\ee
where
\be\label{pbar}
\bar p_t(x):=\frac {\sqrt d} {(2\pi t)^{d/2}}\, e^{-\frac{d|x|^2}{2t}}
\ee
\et
The way we generalize this theorem is by means of the following 
\bt[LCLT for Continuous-Time  Random Walk]\label{LCLTCRW}
Let $x\in \mathbb Z^d$ and  $ p^{\text{RW}}_t(\cdot)$ be the probability distribution of a continuous-time random walk in $\mathbb Z^d$, then, for any fixed $M\ge 0$ there exists $c=c(M)>0$ s.t.
\be
\sup_{|x|\le M \sqrt{t}} \bigg| \frac{p^{RW}_t(x)}{\bar p_t(x)}-1\bigg| \le \frac c {\sqrt t}
\ee
\et

\bpr
We can always decompose 
\be\label{decompose}
p^{\text{RW}}_t(x)=\sum_{n=0}^\infty P(N_t=n) \,  p^{\text{DRW}}_n(x)
\ee
with $N_t$  a Poisson process of rate 1. First by Proposition 2.5.5 in \cite{lawler2010random} we have
\be \label{Poisson1}
P(N_t=n)=\frac 1 {\sqrt{2\pi t}}\, e^{-\frac{(n-t)^2}{2t}} \, \exp\left\{\mathcal{O}\(\frac 1 {\sqrt t}+\frac{|n-t|^3}{t^2}\)\right\} 
\ee
Now for $\epsilon >0$, we assume that  
\be\label{epsassum}
\frac{|n-t|}{t} \leq \epsilon  \nn
\ee
after some manipulation  we obtain the following relations 
\be\label{oneovern}
\frac 1 n = \frac 1 t \(1+\mathcal O\(\frac {|n-t|}{t}\)\), \qquad \frac 1{ n^{\alpha}} = \frac 1 {t^{\alpha}} \(1+\mathcal O\(\frac {|n-t|}{t}\)\)
\ee 
combining \eqref{oneovern} with Theorem \ref{LCLTDRW} we have 
\beq\label{Pndrw}
p^{\text{DRW}}_n(x) &&=\frac {\sqrt d} {(2\pi n)^{d/2}}\, e^{-\frac{d|x|^2}{2n}}  \(1+\mathcal O\(\frac 1 {n}\)\) \nn \\
&&=\frac {\sqrt d} {(2\pi t)^{d/2}}\, e^{-\frac{d|x|^2}{2n}}  \exp \left\{\mathcal O\(\frac{|x|^2|n-t|}{t^2}\)\right\} \(1+\mathcal O\(\frac 1 {t}\)\) \(1+\mathcal O\(\frac {|n-t|}{t}\)\) \nn \\
\eeq
Finally, substitution of \eqref{Poisson1} and \eqref{Pndrw} in \eqref{decompose} and further manipulations gives
\beq
&&\sum_{n=0}^\infty P(N_t=n)  p^{\text{DRW}}_n(x) \nn\\
&&=\sum_{n=0}^\infty \frac 1 {\sqrt{2\pi t}}\, e^{-\frac{(n-t)^2}{2t}} \, \exp\left\{\mathcal O\(\frac 1 {\sqrt t}+\frac{|n-t|^3}{t^2}\)\right\} \; \nn \\
&&\times \frac {\sqrt d} {(2\pi t)^{d/2}}\, e^{-\frac {d|x|^2}{2t}}\, \exp \left\{\mathcal O\(\frac{|x|^2|n-t|}{t^2}\)\right\} \(1+\mathcal O\(\frac 1 {t}\)\) \(1+\mathcal O\(\frac {|n-t|}{t}\)\)\nn \\
\eeq
Assuming $ |x|\le M \sqrt{t}$ and using \eqref{epsassum}, we get the following, 
\be
\exp \left\{\mathcal O\(\frac{|x|^2|n-t|}{t^2}\)\right\} = \exp \left\{\mathcal O\(\epsilon\)\right\}
\ee
Hence, more applications of \eqref{epsassum} give
\beq
&&\sum_{n=0}^\infty P(N_t=n)  p^{\text{DRW}}_n(x) \nn\\
&&= \(1+\mathcal O\(\frac 1 {t}\)\) \frac {\sqrt d} {(2\pi t)^{d/2}}  e^{-\frac {d|x|^2}{2t}} \exp \left\{\mathcal O\(\epsilon\)\right\} \(1+\mathcal O\(\epsilon \)\) \exp\left\{\mathcal O\(\frac 1 {\sqrt t}\)\right\}   \nn \\
&&\times  \sum_{n=0}^\infty \frac 1 {\sqrt{2\pi t}}\, e^{-\frac{(n-t)^2}{2t}} \, \exp\left\{\mathcal O\(\frac{|n-t|^3}{t^2}\)\right\} \; \nn \\
&&= \bar p_t(x)\(1+\mathcal O\(\frac 1 {\sqrt t}\)\)\nn
\eeq
\epr

\section*{Acknowledgements}
 M. Ayala acknowledges financial support from the Mexican Council on Science and Technology (CONACYT)  via the scholarship 457347.

\bibliographystyle{plain}
\bibliography{../Bibliography/Biblio}

\end{document}